\pdfoutput=1
\RequirePackage{ifpdf}
\ifpdf 
\documentclass[pdftex]{sigma}
\else
\documentclass{sigma}
\fi

\usepackage{mathabx}
\usepackage{dcpic,pictexwd}

\newcommand\bfq{\mathbf{q}}
\newcommand\bfp{\mathbf{p}}

\renewcommand\sl{\mathfrak{s}\mathfrak{l} }

\newcommand\CalI{\mathcal{I}}

\newcommand\CalK{\mathcal{K}}

\newcommand\diag{\text{\rm diag}}

\newcommand\pr{\text{\rm pr}}
\newcommand\Real{\mathbb{R}}
\newcommand\CalS{\mathcal{S}}
\DeclareMathOperator{\spn}{span}
\DeclareMathOperator{\rank}{rank}

\newcommand\bfGamma{\boldsymbol{\Gamma}}
\newcommand\ann{\text{\rm ann}}
\newcommand\lag{{{\mathfrak g}}}

\newcommand\homega{\widehat{\omega}}
\newcommand{\hpi}{\widehat{\pi}}

\newcommand\htau{\widehat{\tau}}

\newcommand\cW{\widecheck{W}}

\newcommand\hV{{\widehat{V}}}

\newcommand\hW{\widehat{W}}

\newcommand\ctau{\widecheck{\tau}}

\newcommand\cpi{\widecheck{\pi}}
\newcommand\comega{\widecheck{\omega}}

\newcommand\cV{\widecheck{V}}

\numberwithin{equation}{section}

\newtheorem{Theorem}{Theorem}[section]

\newtheorem{Lemma}[Theorem]{Lemma}
{\theoremstyle{definition}
\newtheorem{Example}[Theorem]{Example}
\newtheorem{Remark}[Theorem]{Remark}
}

\def\beginDC#1[#2]{\begindc{#1}[#2]}

\begin{document}

\allowdisplaybreaks

\renewcommand{\thefootnote}{$\star$}

\renewcommand{\PaperNumber}{017}

\FirstPageHeading

\ShortArticleName{Cauchy Problems and~Darboux Integrability}

\ArticleName{The Cauchy Problem for Darboux Integrable Systems
\\
and Non-Linear d'Alembert Formulas\footnote{This
paper is a~contribution to the Special Issue ``Symmetries of Dif\/ferential Equations: Frames,
Invariants and~Applications''.
The full collection is available
at
\href{http://www.emis.de/journals/SIGMA/SDE2012.html}{http://www.emis.de/journals/SIGMA/SDE2012.html}}}

\Author{Ian M.~ANDERSON and~Mark E.~FELS}

\AuthorNameForHeading{I.M.~Anderson and~M.E.~Fels}

\Address{Utah State University, Logan Utah, USA}

\Email{\href{mailto:Ian.Anderson@usu.edu}{Ian.Anderson@usu.edu},
\href{mailto:Mark.Fels@usu.edu}{Mark.Fels@usu.edu}}

\ArticleDates{Received October 08, 2012, in f\/inal form February 20, 2013; Published online February 27, 2013}

\Abstract{To every Darboux integrable system there is an associated Lie group $G$ which is
a~fundamental invariant of the system and~which we call the Vessiot group.
This article shows that solving the Cauchy problem for a~Darboux integrable partial dif\/ferential
equation can be reduced to solving an equation of Lie type for the Vessiot group~$G$.
If the Vessiot group~$G$ is solvable then the Cauchy problem can be solved by quadratures.
This allows us to give explicit integral formulas, similar to the well known d'Alembert's formula
for the wave equation, to the initial value problem with generic non-characteristic initial data.}


\Classification{58A15; 35L52; 58J70; 35A30; 34A26}
\begin{flushright}
\begin{minipage}{95mm}\it
Dedicated to our good friend and~collaborator Peter Olver on the occasion of his $60^{th}$ birthday.
\end{minipage}
\end{flushright}

\renewcommand{\thefootnote}{\arabic{footnote}}
\setcounter{footnote}{0}

\section{Introduction}

The solution to the classical wave equation $u_{tt}-u_{xx}=0$ with initial data $u(0,x)=a(x)$
and~$u_t(0,x)=b(x)$ is given by the well-known d'Alembert's formula
\begin{gather}
u(t,x)=\frac{1}{2}\big(a(x-t)+a(x+t)\big)+\frac{1}{2}\int_{x-t}^{x+t}b(\xi)d\xi.
\label{DalW}
\end{gather}
In this article we characterize a~broad class of dif\/ferential equations where the solution to the
Cauchy problem can be expressed in terms of the initial data by quadratures as in~\eqref{DalW}.

The family of equations which we identify that can be solved in this manner are a~subset of the
partial dif\/ferential equations known as Darboux integrable equations.
The results we present here are for the classical case of a~scalar Darboux integrable equation in
the plane but we also illustrate how these results do hold in the more general case of a~Darboux
integrable exterior dif\/ferential system.
An example of such an equation is the non-linear hyperbolic PDE in the plane
\begin{gather*}
u_{xy}=\frac{u_x u_y}{u-x}.
\end{gather*}
With initial data given along $y=x$ by $u(x,x)=f(x)$ and~$u_x(x,x)=g(x)$ we f\/ind that the
analogue to~\eqref{DalW} (in null coordinates) is
\begin{gather*}
u=x+\big(f(y)-y\big)e^{\int_x^y G(t)dt}+e^{-\int_0^x G(t)dt}\left(\int_x^y e^{\int_0^s G(t)dt}ds\right),
\end{gather*}
where
\begin{gather*}
G(t)=\frac{g(t)}{t-f(t)}.
\end{gather*}

A fundamental invariant of a~Darboux integrable system, called the Vessiot group, was introduced
in~\cite{anderson-fels-vassiliou:2009a}.
The Vessiot group $G$ is a~Lie group which plays an essential role in the analysis of many of the
geometric properties of Darboux integrable
equations~\cite{anderson-fels:2011a,anderson-fels:2011b}.
This is also true when solving the Cauchy problem for these systems as can be seen by the following
theorem.
\begin{Theorem}
\label{PDEthm}
Let $F(x,y,u,u_x,u_y,u_{xx},u_{xy},u_{yy})=0$ be a~hyperbolic PDE in the plane.
If~$F$ $($or its prolongation$)$ is Darboux integrable then the initial value problem can be solved by
integrating an equation of fundamental Lie type for the Vessiot group~$G$.
If~$G$ is simply connected and~solvable, then the initial value problem can be solved by
quadratures.
\end{Theorem}

Theorem~\ref{PDEthm} can be viewed as a~generalization of the classical theorem of Sophus Lie on
solving an ordinary dif\/ferential equation by quadratures~\cite{olver:1998a}.
This theorem states that the general solution, or the solution to the initial value problem, of an
$n^{\rm th}$ order ordinary dif\/ferential equation with an $n$-dimensional solvable symmetry group can
be solved by quadratures.
While this classical result on integrating ODE's has motivated our work, there is one fundamental
dif\/ference between this classical result for ODE and~Theorem~\ref{PDEthm}.
In Theorem~\ref{PDEthm} the group $G$ is {\it not} a~symmetry group of the PDE $F=0$.

In~\cite{bryant-griffiths-hsu:1995a} it is shown how the initial value problem for Darboux
integrable hyperbolic systems can be solved using the Frobenius theorem.
The approach we take here is quite dif\/ferent.
By using the quotient representation for Darboux integrable hyperbolic Pfaf\/f\/ians systems
constructed in~\cite{anderson-fels-vassiliou:2009a}, we show that the initial value problem can be
solved by solving an equation of fundamental Lie type for the Vessiot group $G$.
This, in turn, allows us to conclude that if the group $G$ is solvable, then the initial value
problem can be solved by quadratures.
We illustrate this with a~number of examples.
The relationship between our approach and~the approach given
in~\cite{bryant-griffiths-hsu:1995a} is described in Appendix~\ref{CCP}.
An example given in the appendix compares the two approaches.

\section{Pfaf\/f\/ian systems and~reduction}\label{section2}

In this section we give the def\/inition of a~Pfaf\/f\/ian system and~summarize some basic facts
about their reduction by a~symmetry group.

\subsection{Pfaf\/f\/ian systems}

A constant rank Pfaf\/f\/ian system is given by a~constant rank sub-bundle $I\subset T^*M$.
An integral manifold of $I$ is a~smooth immersion $s:N\to M$ such that $s^*I=0$.
If $N$ is an open interval, then we call $s$ an integral curve of
$I$,~\cite{bryant-chern-gardner-griffiths-goldschmidt:1991a}.

A {\it local first integral} of a~Pfaf\/f\/ian system $I$ is a~smooth function
$F:U\to\Real$, def\/ined on an open set $U\subset M$, such that $d F\in I$.
For each point $x\in M$ we def\/ine
\begin{gather*}
I_x^\infty=\{dF_x\,|\,\text{$F$ is a local first integral, defined about $x$}\}.
\end{gather*}
We shall assume that $I^\infty={\bigcup}_{x\in M}I_x^\infty$ is a~constant rank bundle on~$M$.
It is easy to verify that $I^\infty$ is the (unique) maximal, completely integrable, Pfaf\/f\/ian
sub-system of~$I$.	Granted additional regularity conditions (see below), the bundle $I^\infty$ can
be computed algorithmically using the derived sequence of~$I$.

The {\it derived system} $I'\subset I$ of a~Pfaf\/f\/ian system $I$ is def\/ined pointwise by
\begin{gather*}
I'_x=\spn\{ \theta_x\,|\,\theta\in\CalS(I) \ {\rm and} \ d\theta\equiv0 \mod I \},
\end{gather*}
where $\CalS(I)$ are the sections of $I$.
The system $I$ is integrable if it satisf\/ies the Frobenius condition $I'=I$.
Letting $I^{(0)}=I$, and~assuming $I^{(k)}$ is constant rank, we def\/ine the derived sequence
inductively by
\begin{gather*}
I^{(k+1)}=\big(I^{(k)}\big)',\qquad k=0,1,\ldots,N,
\end{gather*}
where $N$ is the smallest integer such that $I^{(N+1)}=I^{(N)}$.
Therefore $I^\infty=I^{(N)}$ whenever the sets $I^{(k)}$ are constant rank bundles.
More information about the derived sequence can be found
in~\cite{bryant-chern-gardner-griffiths-goldschmidt:1991a}.

\subsection{Reduction of Pfaf\/f\/ian systems}

A Lie group $G$ acting on $M$ is a~{\it symmetry group} of the Pfaf\/f\/ian system $I$ if
\begin{gather*}
g^*I=I, \qquad{\rm for~all} \ g\in G.
\end{gather*}
The group $G$ acts \emph{regularly} on $M$ if the quotient map $\bfq_G:M\to M/G$ is a~smooth submersion.
Let $\Gamma$ denote the Lie algebra of inf\/initesimal generators for the action of~$G$ on~$M$
and~let $\bfGamma\subset TM$ be its pointwise span.
Then $\bfGamma=\ker(\bfq_{G*})$.
Assume from now on that the action of $G$ on $M$ is regular and~a~symmetry group of~$I$.

The reduction $I/G$ of $I$ by the symmetry group $G$ is def\/ined by
\begin{gather}
I/G=\big\{ \bar\theta \in\Lambda^1(M/G)\,|\,\bfq_G^*\bar\theta\in I\big\}.
\label{Iquot}
\end{gather}
Let $\ann(I)\subset TM$ be the annihilating space of $I$,
\begin{gather*}
\ann(I)=\big\{X\in TM\,|\,\theta(X)=0 \ {\rm for~all} \ \theta\in I\big\}.
\end{gather*}
It is easy to check that $I/G$ is constant rank if and~only if $\bfGamma\cap\ann(I)$ is constant
rank~\cite{anderson-fels:2005a} in which case
\begin{gather*}
\rank(I/G)=\rank(I)-\big(\rank(\bfGamma)-\rank(\bfGamma\cap\ann(I))\big)=
\rank(\ann(\bfGamma)\cap I).
\end{gather*}
The annihilating space $\ann(I)$ is $G$ invariant and~satisf\/ies~\cite{fels:2008a},
\begin{gather}
\bfq_{G*}(\ann(I))=\ann(I/G).
\label{redannI}
\end{gather}
For examples on computing $I/G$, applications, and~more information about its properties
see~\cite{anderson-fels:2011a,anderson-fels:2005a,anderson-fels:2011b,fels:2008a}.

The symmetry group $G$ is said to act {\it transversally} to the Pfaf\/f\/ian system $I$ if
\begin{gather}
\ker(\bfq_{G*})\cap\ann(I)=0.
\label{TransC2}
\end{gather}
The following basic theorem, which follows directly from Theorem~2.2 in~\cite{anderson-fels:2011a},
illustrates the importance of the transversality condition~\eqref{TransC2}.
\begin{Theorem}
\label{FRD}
Let $G$ be a~symmetry group of the constant rank Pfaffian system $I$ on $M$ acting freely
and~regularly on $M$ and~transversally to $I'$.
Then $I/G$ is a~constant rank Pfaffian system and~$(I')/G=(I/G)'$.
Furthermore, if $\gamma:(a,b)\to M/G$ is a~one-dimensional integral manifold of $I/G$,
then through each point $p\in M$ satisfying $\bfq_G(p)=\gamma(t_0)$ there exists a~unique lift
$\sigma:(a,b)\to M$ of $\gamma$ which is an integral manifold of $I$ satisfying $\sigma(t_0)=p$.
The lifted curve $\sigma$ can be found by solving an equation of fundamental Lie type.
\end{Theorem}

Let $\lag$ be the Lie algebra of $G$.
Then, given a~curve $\alpha:(a,b)\to\lag$, the system of ODE
\begin{gather}
\dot\lambda(t)=(L_\lambda)_*(\alpha(t))
\label{LT}
\end{gather}
for the curve $\lambda:\Real\to G$ is called an {\it equation of fundamental Lie type}.
Given $\alpha$ in equation~\eqref{LT}, the curve $\lambda$ can be found globally over the entire
interval $(a,b)$.
See~\cite[p.~55]{bryant:1991a} for more information about equations
of Lie type.
If $G$ is simply connected and~$\lag$ is solvable, then an equation of fundamental Lie type can be
solved by quadratures, see~\cite[Proposition~4, p.~60]{bryant:1991a}.
\begin{proof}
The f\/irst statement in Theorem~\ref{FRD} is proved
in~\cite{anderson-fels:2005a}.
For the second statement concerning the lifting of integral curves $\gamma$ we begin by noting that
since $G$ is free we may, by using a~$G$-invariant Riemannian metric on $M$, write
\begin{gather}
I=\bfq_G^*(I/G)\oplus K,
\label{Iplus}
\end{gather}
where $\bfq^*(I/G)$ is the pullback bundle and~$K$ is its $G$-invariant orthogonal complement in
$I$.
The transversality condition~\eqref{TransC2} implies that
\begin{gather*}
\rank(K)=\dim G\qquad{\rm and}\qquad\ann(K)\cap\ker(\bfq_{G*})=0.
\end{gather*}
Therefore $\ann(K)$ is a~horizontal space for the principal bundle $M\to M/G$.
Consequently, since $\gamma$ is an integral curve of $I/G$, a~lift $\sigma$ which is an integral
curve of $I$ satisf\/ies on account of equation~\eqref{Iplus}, $\sigma^*(K)=0$.
Therefore $\sigma$ is a~horizontal lift for the connection $\ann(K)$ on the principal bundle $M$.
But the system of ODE determining a~horizontal lift for a~connection on a~principal bundle is
precisely an equation of fundamental Lie type for the structure group~$G$, see~\cite[p.~69]{kobayashi-nomizu:2009a}.
\end{proof}
\begin{Remark}
\label{solvele}
There are two standard ways to f\/ind a~lift $\sigma:(a,b)\to M$ of $\gamma$.
To describe them let $p\in M$ satisfy $\bfq_{G}(p)=\gamma(t_0)$.
We now describe how to f\/ind the lift satisfying $\sigma(t_0)=p$.

The f\/irst way is to start with {\it any} lift $\hat\sigma:(a,b)\to M$ of $\gamma$ satisfying
$\hat\sigma(t_0)=p$.
Then def\/ine a~lift $\sigma(t)=\lambda(t)\hat\sigma(t)$, where $\lambda:(a,b)\to G$
and~$\lambda(t_0)=e_G$.
The requirement that $\sigma$ be an integral curve of $I$ is equivalent to $\lambda$ satisfying an
equation of fundamental Lie type.
The equation of Lie type is easily determined algebraically from $I$ and~$\hat\sigma$.

The second way to f\/ind the lift $\sigma$ is to assume $\gamma:(a,b)\to M/G$ is an embedding.
Then let $S=\gamma(a,b)$ and~$P=\bfq_G^{-1}(S)$.
Theorem~\ref{FRD} states that the restriction $(I)_P$ is the Pfaf\/f\/ian system for an equation
of fundamental Lie type in the sense of~\cite{doubrov:2000a} or~\cite{lych:1991a}.
We then compute the maximal integral manifold through $p\in P$ in order to f\/ind $\gamma$.
\end{Remark}

\section{Hyperbolic Pfaf\/f\/ian systems}\label{section3}

A constant rank Pfaf\/f\/ian system $I$ is said to be hyperbolic of class $s$
(see~\cite{bryant-griffiths-hsu:1995a}) if the following holds.
About each point $x\in M$ there exists an open set $U$ and~a~coframe
$\{ \theta^1,\ldots,\theta^s,
\homega, \hpi, \comega, \cpi \}$ on $U$
such that
\begin{gather*}
I=\spn\big\{ \theta^1,\ldots,\theta^s \big\}
\end{gather*}
and the following structure equations hold,
\begin{alignat*}{3}
& d\theta^{i}\equiv0 \quad && \mod I,\qquad1\leq i\leq s-2, &
\\
& d\theta^{s-1}\equiv\homega\wedge\hpi\quad  && \mod I, &
\\
& d\theta^{s}\equiv\comega\wedge\cpi \quad && \mod I.&
\end{alignat*}

The two Pfaf\/f\/ian systems $\hV$ and~$\cV$ def\/ined by
\begin{gather}
\hV=\spn\big\{ \theta^i, \hpi, \homega \big\},\qquad
\cV=\spn\big\{ \theta^i, \cpi, \comega\big\}
\label{SPS}
\end{gather}
are called the associated {\it characteristic or singular systems} of~$I$.
They are important invariants of~$I$.
Let $\hV^\infty$ and~$\cV^\infty$ be the corresponding subspaces of f\/irst integrals for the
singular Pfaf\/f\/ian systems $\hV$ and~$\cV$ in~\eqref{SPS} of the hyperbolic Pfaf\/f\/ian system~$I$.
A locally def\/ined function $F$ satisfying $dF\in\hV^\infty$ or $dF\in\cV^\infty$ is called an
intermediate integral of $I$ or a~Darboux invariant of~$I$.
We will call $\hV^\infty$ and~$\cV^\infty$ the spaces of intermediate integrals of~$I$.
The hyperbolic system~$I$ is said to be {\it Darboux integrable} if $I^\infty=0$, and
\begin{gather}
\hV+\cV^\infty=T^*M \qquad{\rm and}\qquad\hV^\infty+\cV=T^*M.
\label{DIC}
\end{gather}
See~\cite{anderson-fels-vassiliou:2009a} and~Theorem~4.3 in~\cite{anderson-fels:2011a}.
In particular
Theorem~4.3 in~\cite{anderson-fels:2011a} shows that the conditions in~\eqref{DIC} imply
\begin{gather*}
\hV^\infty\cap\cV^\infty=0.
\end{gather*}

For more information about hyperbolic Pfaf\/f\/ian systems see the
article~\cite{bryant-griffiths-hsu:1995a} where the ge\-ne\-ral theory of hyperbolic exterior
dif\/ferential systems is developed in detail.
The def\/inition of a~hyperbolic exterior dif\/ferential system is generalized
in~\cite{anderson-fels-vassiliou:2009a}, see Section~\ref{decomp}.

The characteristic directions for a~hyperbolic system $I$ with singular systems $\hV$ and~$\cV$ are
$\ann\big(\hV\big)$ and~$\ann(\cV)$~\cite{bryant-griffiths-hsu:1995a}.
A {\it non-characteristic integral curve} is an immersion $\gamma:(a,b)\to M$ which is
a~one-dimensional integral manifold of $I$ such that
\begin{gather}
\dot\gamma(t)\not\in\ann\big(\hV\big)\qquad{\rm and}\qquad\dot\gamma(t)\not\in\ann(\cV),\qquad{\rm for~all} \ t\in(a,b).
\label{nonchar}
\end{gather}
Given a~non-characteristic integral curve $\gamma:(a,b)\to M$ of $I$, a~solution to the Cauchy or
initial value problem for $\gamma$ is a~$2$-dimensional integral manifold $s:N\to M$ of $I$ such
that $\gamma(a,b)\subset s(N)$.
A local solution to the Cauchy problem about a~point $x\in\gamma(a,b)$ is an open neighbourhood
$U\subset M$, $x\in U$, and~a~$2$-dimensional integral manifold $s:N\to U$ of~$I$ such that
$\gamma(a,b)\cap U\subset s(N)$.

Hyperbolic Pfaf\/f\/ian systems of class $s=3$ are closely related to hyperbolic partial
dif\/ferential equations in the plane~\cite{gardner-kamran:1993a}.
Specif\/ically, a~hyperbolic PDE in the plane
\begin{gather}
F(x,y,u,u_x,u_y,u_{xx},u_{xy},u_{yy})=0
\label{FPDE}
\end{gather}
def\/ines a~7-dimensional submanifold $M\subset J^2(\Real^2,\Real)$.
Let $I$ be the rank 3 Pfaf\/f\/ian system which is the restriction of the rank 3 contact system on
$J^2(\Real^2,\Real)$ to $M$.
The following structure equations proven in~\cite{gardner-kamran:1993a} show that $I$ is a~class
$s=3$ hyperbolic Pfaf\/f\/ian system.
\begin{Theorem}
\label{hyperpde}
Let $I$ be the rank $3$ Pfaffian system on the $7$-manifold~$M$ defined by the hyperbolic PDE
in the plane in~\eqref{FPDE}.
About each point $p\in M$ there exists a~local coframe $\{\theta^i,\omega^a,\pi^a\}_{0\leq i\leq2,1\leq a\leq2}$ such that
\begin{gather*}
I=\spn\big\{ \theta^0, \theta^1, \theta^2 \big\},
\end{gather*}
and
\begin{alignat}{3}
& d\theta^0=\theta^1\wedge\omega^1+\theta^2\wedge\omega^2\quad && \mod\big\{\theta^0\big\}, & \nonumber
\\
& d\theta^1=\omega^1\wedge\pi^1+\mu_1\theta^2\wedge\pi^2\quad && \mod\big\{\theta^0,\theta^1\big\},& \nonumber
\\
& d\theta^2=\omega^2\wedge\pi^2+\mu_2\theta^1\wedge\pi^1\quad && \mod\big\{\theta^0,\theta^2\big\}, &
\label{NKstruct}
\end{alignat}
where $\mu_1$, $\mu_2$ are the Monge--Ampere invariants.
The invariant conditions $\mu_1=\mu_2=0$ are satisfied if and~only if~\eqref{FPDE} is locally
a~Monge--Ampere equation.
\end{Theorem}

The article by Gardner and~Kamran~\cite{gardner-kamran:1993a} gives an algorithm to construct the
coframe in equations~\eqref{NKstruct}.
The corresponding singular systems $\hV$ and~$\cV$ are denoted by $C(\CalI_F,dM_1)$
and~$C(\CalI_F,dM_2)$ in~\cite{gardner-kamran:1993a}, and~in frame~\eqref{NKstruct},            
$\hV=\{\theta^0,\theta^1,\theta^2,\omega^1,\pi^1\}$
and~$\cV=\{\theta^0,\theta^1,\theta^2,\omega^2,\pi^2\}$.

\section{The Vessiot group and~the quotient representation}\label{section4}

\subsection{The quotient representation}

We now explain the fundamental role played by the quotient
of a~Pfaf\/f\/ian system by a~symmetry group in the theory of Darboux integrability.
See~\cite[Section~6, p.~20]{anderson-fels:2011a} (with $L=G_{\diag}$) for more information
and~detailed proofs of the following key theorem.
\begin{Theorem}
\label{Dreduce}
Let $K_i$, $i=1,2$ be constant rank Pfaffian systems on $M_i$, $i=1,2$ of codimension~$2$,
and~satisfying $K_i^\infty=0$.
Consider a~Lie group $G$ which acts freely and~regularly on $M_i$, is a~{\it common} symmetry
group of both $K_1$ and~$K_2$ and~acts transversely to $K_1$ and~$K_2$.
Assume also that the action of the diagonal subgroup $G_{\diag}\subset G\times G$ on $M_1\times
M_2$ is regular and~acts transversely to $K'_1+K'_2$.
\begin{itemize}\itemsep=0pt
\item[$(i)$] The sum $K_1+K_2$ on $M_1\times M_2$ is a~constant rank, Darboux integrable, hyperbolic
Pfaffian system.
\item[$(ii)$] The singular Pfaffian systems for $K_1+K_2$ are
\begin{gather*}
\hspace*{-0.9cm}
\hW=K_1+T^*M_2\qquad  and \qquad\cW=T^*M_1+K_2.
\end{gather*}
\item[$(iii)$] The quotient differential system $I=(K_1+K_2)/G_{\diag}$ on $M=(M_1\times
M_2)/G_{\diag}$	is a~constant rank hyperbolic Pfaffian system which is Darboux integrable.
\item[$(iv)$] The singular Pfaffian systems for the quotient system $I$ are
\begin{gather*}
\hspace*{-0.9cm}
\hV=\big(K_1+T^*M_2\big)/G_\diag=\hW/G_\diag\qquad{and}\qquad
\cV=\big(T^*M_1+K_2\big)/G_{\diag}=\cW/G_\diag.
\end{gather*}
\item[$(v)$] The spaces of intermediate integrals for $I$ are
\begin{gather*}
\hspace*{-0.9cm}
\hV^\infty=(0+T^*M_2)/G_\diag=\hW^\infty/G_\diag\!\qquad{and}\!\qquad
\cV^\infty=(T^*M_1+0)/G_{\diag}=\cW^\infty/G_\diag.
\end{gather*}
\end{itemize}
\end{Theorem}

The sums of the type $K_1+K_2$ in Theorem~\ref{Dreduce} are def\/ined precisely as
\begin{gather}
K_1+K_2=\pi_1^*(K_1)+\pi_2^*(K_2),
\label{K1pK2}
\end{gather}
where $\pi_i:M_1\times M_2\to M_i$.

Theorem~\ref{Dreduce} shows how Darboux integrable hyperbolic Pfaf\/f\/ian systems can be
constructed using the group reduction of pairs of Pfaf\/f\/ian systems.
It is a~remarkable fact, established in~\cite{anderson-fels-vassiliou:2009a}, that the converse is
true locally, that is, every Darboux integrable hyperbolic system can be realized locally as
a~non-trivial quotient of a~pair of Pfaf\/f\/ian systems with a~common symmetry group.
The precise formulation of this result is as follows.
\begin{Theorem}
\label{Intro4}
Let $I$ be a~Darboux integrable hyperbolic Pfaffian system on a~manifold $M$ and let~$\hV$
and~$\cV$ be the singular Pfaffian systems as in~\eqref{SPS}.
Fix a~point $x_0$ in $M$ and~let
\begin{itemize}\itemsep=0pt
\item[$(i)$] $M_1$ and~$M_2$ be the maximal integral manifolds of $\hV^\infty$ and~$\cV^\infty$
through $x_0$, and
\item[$(ii)$] $K_1$ and~$K_2$ be the restrictions of $\hV$ and~$\cV$ to $M_1$ and~$M_2$ respectively.
\end{itemize}
\noindent
Then there are open sets $U\subset M$, $U_1\subset M_1$, $U_2\subset M_2$, each containing $x_0$,
and a~local action of a~Lie group $G$ on $U_1$ and~$U_2$ which satisfy the hypothesis of
Theorem~{\rm \ref{Dreduce}} and~such that
\begin{gather}
U=(U_1\times U_2)/G_\diag,\qquad
I_U=\left((K_1+K_2)_{U_1\times U_2}\right)/G_\diag,
\label{LocQuotRep}
\end{gather}
and properties $(iv)$ and~$(v)$ in Theorem~{\rm \ref{Dreduce}} hold on $U$.
\end{Theorem}

The group $G$ appearing in Theorems~\ref{Dreduce} and~\ref{Intro4} is called the {\it
Vessiot group} of the Darboux integrable system $I$.
We shall refer to~\eqref{LocQuotRep} or~$(iii)$ in Theorem~\ref{Dreduce} as {\it the canonical
quotient representation for a~Darboux integrable hyperbolic Pfaffian system~$I$}.
\begin{Remark} It is an algorithmic (but non-trivial) process to f\/ind the group $G$ and~its
action in Theorem~\ref{Intro4}.
The Lie algebra of inf\/initesimal generators of $G$ can be found algebraically, while f\/inding
the action of $G$ may require solving a~system of Lie
type~\cite{carinena:1995a,doubrov:2000a}.
\end{Remark}

\section[Solving the Cauchy initial value problem for Darboux integrable systems]{Solving the Cauchy initial value problem\\ for Darboux integrable systems} \label{section5}

In this section we solve the initial value problem for a~Darboux integrable hyperbolic Pfaf\/f\/ian
system and~give a~proof of Theorem~\ref{PDEthm}.
In the second subsection we outline how to extend these results to the general case of Darboux
integrable systems which are not necessarily hyperbolic.

\subsection{The Cauchy problem for hyperbolic systems}
\label{CPS}

We begin with a~key lemma which shows that the lift of a~non-characteristic integral curve is again
a~non-characteristic integral curve.
\begin{Lemma}
\label{liftCC}
Let $I$ be a~hyperbolic Pfaffian system which is Darboux integrable and~let
$(K_1+K_2)/G_{\diag}$ be the canonical quotient representation of $I$ as in Theorem~{\rm \ref{Intro4}}.
Let $\gamma:(a,b)\to M$ be a~non-characteristic integral curve of $I$ and~let $\sigma:(a,b)\to
M_1\times M_2$ be a
lift of~$\gamma$ to a~one-dimensional integral curve of $K_1+K_2$.
Then~$\sigma$ is a~non-characteristic integral curve for $K_1+K_2$.
\end{Lemma}

The existence of the lifted integral curve $\sigma$ in Lemma~\ref{liftCC} is guaranteed by
Theorem~\ref{FRD}.
\begin{proof} In parts $(i)$ and~$(ii)$ of Theorem~\ref{Dreduce} it is noted that $K_1+K_2$
is Darboux integrable with characteristic systems $\hW=(K_1+T^*M_2)$ and~$\cW=T^*M_1+K_2$.
Therefore, in view of equation~\eqref{nonchar}, we prove Lemma~\ref{liftCC} by showing
\begin{gather*}
\dot\sigma(t)\not\in\ann(\hW)\qquad{\rm and}\qquad\dot\sigma(t)\not\in\ann(\cW),\qquad{\rm for~all} \ t\in(a,b).
\end{gather*}

From part $(iv)$ of Theorem~\ref{Dreduce} and~equation~\eqref{redannI} we f\/ind that
\begin{gather}
\bfq_{G_\diag*}\big(\ann(\hW)\big) =\ann\big(\hW/G_\diag\big)=\ann\big(\hV\big)\qquad{\rm and}
\nonumber\\
 \bfq_{G_\diag*}\big(\ann(\cW)\big) =\ann\big(\cW/G_\diag\big)=\ann\big(\cV\big).
\label{PFV}
\end{gather}
Now suppose that $\sigma$ is characteristic at some point $t_0$ so that, for instance,
$\dot\sigma(t_0)\in\ann(\hW)$.
Then by~\eqref{PFV}
\begin{gather*}
\bfq_{G_\diag*}\big(\dot\sigma(t_0)\big)\in\ann\big(\hV\big)
\end{gather*}
and hence, since $\sigma$ is a~lift of $\gamma$,
$\dot\gamma(t_0)=\bfq_{G_\diag*}\big(\dot\sigma(t_0)\big)\in\ann\big(\hV\big)$.
This contradicts the fact that~$\gamma$ is not characteristic.
A~similar argument applies if we assume $\dot\sigma(t_0)\in\cW$.
\end{proof}

We also need the following lemma.
\begin{Lemma}
\label{fundl}
Let $\sigma:(a,b)\to M_1\times M_2$ be a~non-characteristic integral curve of $K_1+K_2$, where~$K_1$ and~$K_2$ satisfy the conditions of Theorem~{\rm \ref{Dreduce}}.
Then the curves $\sigma_i:(a,b)\to M_i$ defined by
\begin{gather}
\sigma_i=\pi_i\circ\sigma
\label{sigmai}
\end{gather}
are $1$-dimensional integral manifolds of $K_i$.
\end{Lemma}

\begin{proof} By applying~\eqref{K1pK2} and~\eqref{sigmai} we f\/ind
\begin{gather*}
\sigma_i^*(K_i)=\sigma^*\pi_i^*(K_i)\subset\sigma^*(K_1+K_2)=0.
\end{gather*}
It then remains to be shown that the maps $\sigma_i:(a,b)\to M_i$ are immersions.

Suppose $\dot\sigma_1(t_0)=\left(\pi_{1*}(\dot\sigma)\right)(t_0)=0$.
Then on the one-hand,
\begin{gather}
\dot\sigma(t_0)\in\ker(\pi_{1*})=0+TM_2.
\label{cap1}
\end{gather}
On the other hand, by writing
\begin{gather*}
K_1+K_2=\big(K_1+T^*M_2\big)\cap(T^*M_1+K_2)
\end{gather*}
we also have that
\begin{gather}
\dot\sigma(t_0)\in\ann(K_1+K_2)=\ann\big(K_1+T^*M_2\big)\oplus\ann(T^*M_1+K_2).
\label{PL2}
\end{gather}
Therefore, since
\begin{gather*}
\ann\big(K_1+T^*M_2\big)\cap(0+TM_2)=0\qquad{\rm and}\qquad
\ann\big(T^*M_1+K_2\big)\subset(0+TM_2),
\end{gather*}
we get from equations~\eqref{cap1} and~\eqref{PL2}
\begin{gather*}
\dot\sigma(t_0)\in\ann\big(T^*M_1+K_2\big)=\ann(\cW),
\end{gather*}
which contradicts the hypothesis that $\sigma$ is non-characteristic.
A similar argument applies to $\sigma_2$ and~so we conclude that the maps $\sigma_i:(a,b)\to M_i$
are immersions.
\end{proof}

The solution to the Cauchy initial value problem can now be given.
\begin{Theorem}
\label{slnthm}
Let $\gamma:(a,b)\to M$ be a~non-characteristic integral curve for the Darboux integrable
hyperbolic Pfaffian system $I$.
Let $I$ have canonical quotient representation $(K_1+K_2)/G_\diag$, let $\sigma:(a,b)\to M_1\times
M_2$ be a~lift of $\gamma$ to an integral curve of $K_1+K_2$ and~let
$\sigma_i=\pi_i\circ\sigma:(a,b)\to M_i$. Then the smooth function $s:(a,b)\times(a,b)\to M$
defined by
\begin{gather}
s(t_1,t_2)=\bfq_{G_\diag}\big(\sigma_1(t_1),\sigma_2(t_2)\big)
\label{defSig}
\end{gather}
solves the Cauchy problem for $\gamma$.
\end{Theorem}

\begin{proof} We f\/irst show that $s$ is an integral manifold of $I$.
Def\/ine $\Sigma:(a,b)\times(a,b)\to M_1\times M_2$~by
\begin{gather*}
\Sigma(t_1,t_2)=\big(\sigma_1(t_1),\sigma_2(t_2)\big).
\end{gather*}
Lemma~\ref{fundl} implies that $\Sigma$ is a~2-dimensional integral manifold of $K_1+K_2$.
By the def\/inition of quotient in~\eqref{Iquot}, $\bfq_{G_\diag}$ maps integral manifolds to
(possibly non-immersed) integral manifolds.
The condition that $G_{\diag}$ acts transversally to $K_1+K_2$ (condition~\eqref{TransC2})
guarantees that since~$\Sigma$ is an integral of $K_1+K_2$, the composition
$\bfq_{G_\diag}\circ\Sigma$ is an immersion.
Therefore $s=\bfq_{G_\diag}\circ\Sigma$ is a~2-dimensional integral manifold of $I$.

We now show that $\gamma(a,b)\subset s((a,b)\times(a,b))$.
All we need to do is set $t_1=t_2=t$ in equation~\eqref{defSig}.
Since $\sigma(t)=(\sigma_1(t),\sigma_2(t))$ we have
\begin{gather*}
\Sigma(t,t)=\bfq_{G_\diag}(\sigma_1(t),\sigma_2(t))=\bfq_{G_\diag}\circ\sigma(t)=\gamma(t)
\end{gather*}
and therefore $\gamma(a,b)\subset s((a,b)\times(a,b))$.
\end{proof}

The proof of Theorem~\ref{PDEthm} is now simple.

\begin{proof} If $F(x,y,u,u_x,u_y,u_{xx},u_{xy},u_{yy})=0$ is a~hyperbolic PDE which is Darboux
integrable after $k$ prolongations then $I^{\langle k\rangle}$, the $k^{\rm th}$ prolongation of the rank 3
Pfaf\/f\/ian system $I$ in Theorem~\ref{hyperpde}, is a~hyperbolic Pfaf\/f\/ian system of class
$3+k$ which is Darboux integrable~\cite{bryant-griffiths-hsu:1995a}.
Non-characteristic initial data $\gamma:(a,b)\to M$ prolongs to non-characteristic initial data
$\gamma^{\langle k\rangle }$ for $I^{\langle k\rangle }$.
We then apply Theorem~\ref{slnthm} to $I^{\langle k\rangle}$ to solve the initial value problem for $I^{\langle k\rangle }$
using the prolonged Cauchy data $\gamma^{\langle k\rangle }$.
The solution~$s$ to the Cauchy problem for the prolongation projects to the solution to the initial
value problem for~$I$ and~hence $F$.

Finally, the solution given in Theorem~\ref{slnthm} only requires computing the lift of the initial
data~$\gamma$ to an integral curve of $K_1+K_2$ on the product space $M_1+M_2$.
This, by Theorem~\ref{FRD}, only involves solving an equation of fundamental Lie type.
\end{proof}

\begin{Remark} Given a~Darboux integrable hyperbolic system $I$, Theorem~\ref{Intro4} only
guarantees that~$I$ admits a~local quotient representation in the sense that the action of~$G$ is
local and~the quotient representation~\eqref{LocQuotRep} only holds locally.
In this case the implementation of Lemmas~\ref{liftCC},~\ref{fundl} and~Theorem~\ref{slnthm}
will only produce a~local solution to the Cauchy problem.
However it is still the case that if the Vessiot group~$G$ is solvable, then this local solution
can be found by quadratures.
\end{Remark}

At this point the reader may wish to refer to Appendix~\ref{CCP} to compare the results in this
subsection with the classical approach to solving the Cauchy problem for Darboux integrable
hyperbolic systems.

\subsection{Generalization to Darboux integrable systems}
\label{decomp}

In this section we outline how to solve the initial value problem in the general case of Darboux
integrable Pfaf\/f\/ian systems which are not necessarily hyperbolic.
A simple demonstration is given in Example~\ref{lastE}.

A Pfaf\/f\/ian system $I$ is a~called decomposable if about each point $p\in M$ there exists an
open set $U\subset M$ and~a~coframe on $U$ given by
\begin{gather*}
\theta^1,\ldots,\theta^r, \
\homega^1,\dots,\homega^{n_1}, \ \htau^1,\dots,\htau^{p_1}, \
\comega^1,\dots,\comega^{n_2}, \ \ctau^1,\dots,\ctau^{p_2},
\end{gather*}
where $n_1+p_1\geq2$, $n_2+p_2\geq2$, $n_1,n_2,p_1,p_2,\geq1$, with the properties (see Theorem 2.3
in~\cite{anderson-fels-vassiliou:2009a}),
\begin{itemize}\itemsep=0pt
\item[$(i)$] the Pfaf\/f\/ian system is $I_U=\spn\{\theta^i\},\;1\leq i\leq r$;
\item[$(ii)$] the structure equations are
\begin{alignat*}{3}
& d\theta^{i_0} \equiv0 \quad && \mod I,\qquad1\leq i_0\leq r_1, &
\\
& d\theta^{i_1} \equiv A^{i_1}_{ab}\htau^a\wedge\homega^b \quad &&  \mod I,\qquad r_1+1\leq i_1\leq r_2,&
\\
& d\theta^{i_2} \equiv
B^{i_2}_{\alpha\beta}\ctau^\alpha\wedge\comega^\beta \quad && \mod I,\qquad r_2+1\leq i_2\leq r;&
\end{alignat*}
\item[$(iii)$] and~$I'=\spn\{\theta^{i_0}\}$.
\end{itemize}

Decomposable systems are generalizations of the hyperbolic systems def\/ined in Section~\ref{section3}.

{\samepage
The Pfaf\/f\/ian systems $\hV$, $\cV$ def\/ined by
\begin{gather*}
\hV=\spn\big\{ \theta^i, \htau^a, \homega^b\big\}\qquad{\rm and}\qquad
\cV=\spn\big\{ \theta^i, \ctau^\alpha, \comega^\beta\big\}
\end{gather*}
are again called the characteristic or singular systems of a~decomposable Pfaf\/f\/ian system~$I$.
As in the case of hyperbolic systems, a~decomposable Pfaf\/f\/ian system $I$ is said to be Darboux
integrable if $I^\infty=0$, and
\begin{gather*}
\hV+\cV^\infty=T^*M,\qquad\hV^\infty+\cV=T^*M.
\end{gather*}
See~\cite{anderson-fels-vassiliou:2009a} and~Theorem~4.3 in~\cite{anderson-fels:2011a}.}

Theorems~\ref{Dreduce} and~\ref{Intro4} hold for Darboux integrable systems by simply dropping the
codimension 2 requirement in Theorem~\ref{Dreduce}, see~\cite{anderson-fels-vassiliou:2009a}
and~\cite{anderson-fels:2011a}.
This provides the canonical quotient representation of $I$ used in the analogue of
Theorem~\ref{liftCC}.

The Cauchy initial value problem for a~decomposable Pfaf\/f\/ian system consists of f\/irst
prescribing an $(n_1+n_2-1)$-dimensional non-characteristic integral manifold $\gamma:S\to M$ which
we will assume to be embedded.
We then need to f\/ind an $(n_1+n_2)$-dimensional integral manifold $s:N\to M$ such that
$\gamma(S)\subset s(N)$.

The full generalization of Theorem~\ref{PDEthm} is the following.
\begin{Theorem}
\label{DIG}
Let $I$ be a~decomposable Darboux integrable Pfaffian system.
The initial value problem for $I$ can be solved by integrating an equation of fundamental Lie type
for the Vessiot group $G$.
If $G$ is simply connected and~solvable, then the initial value problem can be solved by
quadratures.
\end{Theorem}

The steps needed to prove Theorem~\ref{DIG} and~solve the Cauchy problem are identical to those in
Section~\ref{CPS}.
Lemmas~\ref{liftCC},~\ref{fundl} and~Theorem~\ref{slnthm} have direct generalizations which appear
in this section.
However the proofs of these corresponding results are now signif\/icantly dif\/ferent.
This is due to the fact that the condition that $S$ be non-characteristic for a~general
decomposable system is considerably more complex than condition~\eqref{nonchar} for a~curve to be
non-characteristic for a~hyperbolic system.
\begin{Lemma}
\label{liftCC2}
Let $I$ be a~decomposable Darboux integrable Pfaffian system and~$($as in Theorem~{\rm \ref{Intro4})}
let $((K_1+K_2)/G_{\diag},(M_1\times M_2)/G_\diag)$ be the canonical quotient representation of $I$
where $G$ is the Vessiot group.
Let $S\subset M$ be a~non-characteristic integral manifold of $I$ of dimension $n_1+n_2-1$.
Choose $x_0\in S$ and~$(x_1,x_2)\in M_1\times M_2$ with $\bfq_{G_\diag}(x_1,x_2)=x_0$, and let~$L$
be the maximal integral manifold of the fundamental Lie system
\begin{gather*}
(K_1+K_2)|_{\bfq^{-1}_{G_\diag}(S)}
\end{gather*}
through $(x_1,x_2)$.
Then $L\subset M_1\times M_2$ is an $(n_1+n_2-1)$-dimensional non-characteristic integral manifold of $K_1+K_2$.
\end{Lemma}

See~\cite{carinena:1995a} and~\cite{doubrov:2000a} for more information on systems of equations of
Lie type.
The generalization of Lemma~\ref{fundl} is then the following.
\begin{Lemma}
\label{fundl2}
Let $L\subset M_1\times M_2$ be an $(n_1+n_2-1)$-dimensional embedded non-characteristic integral
manifold of $K_1+K_2$, where $K_1$ and~$K_2$ satisfy the conditions of Theorem~{\rm \ref{Dreduce}}.
Then the manifolds $L_i\subset M_i$ defined by
\begin{gather*}
L_i=\pi_i(L)
\end{gather*}
satisfy $\dim L_1=n_2$ and~$\dim L_2=n_1$ and~are integral manifolds of $K_i$.
\end{Lemma}

{\samepage
Applying Lemmas~\ref{liftCC2} and~\ref{fundl2} produces in a~manner similar to
Theorem~\ref{slnthm} the following solution to the initial value problem.
\begin{Theorem}
\label{slnthm3}
Let $S\subset M$ be a~non-characteristic $(n_1+n_2-1)$-dimensional integral manifold for the Darboux
integrable Pfaffian system $I$.
Let $L\subset M_1\times M_2$ be a~lift of $S$ to an integral manifold of $K_1+K_2$ satisfying the
conditions in Lemma~{\rm \ref{liftCC2}}, and~let $L_i=\pi_i(L)$.
Then the smooth function $s:L_1\times L_2\to M$ defined by
\begin{gather*}
s(t_1,t_2)=\bfq_{G_\diag}(t_1,t_2),\qquad t_1\in L_1,\qquad t_2\in L_2,
\end{gather*}
solves the local Cauchy problem for $I$ about $x_0$ with Cauchy data $S\subset M$.
\end{Theorem}

Lemma~\ref{liftCC2} and~Theorem~\ref{slnthm3} establish Theorem~\ref{DIG}.}

\begin{Remark} As this paper was near completion we obtained the preprint~\cite{vassiliou:2012a}
which gives a~detailed example of a~Darboux integrable wave map system where the solution to the
Cauchy problem reduces to the integration of a~Lie system for ${\rm SL}(2,\Real)$.
This is an excellent example of Theorem~\ref{slnthm3}.
\end{Remark}

\section{Examples}
\begin{Example}
For our f\/irst example we consider the Darboux integrable partial dif\/ferential equation
\begin{gather}
u_{xy}=\frac{u_xu_y}{u-x}.
\label{pde1}
\end{gather}
The closed-form general solution to this equation is well known (for example
see~\cite{anderson-fels-vassiliou:2009a}).
Here however we will f\/ind the solution to~\eqref{pde1} with initial Cauchy data by demonstrating
the constructions in Theorems~\ref{Intro4} and~\ref{slnthm}.
Let $\gamma:\Real\to M$ be the Cauchy data given by
\begin{gather}
\gamma(x)=\bigg( x=x, \; y=x, \; u=f(x), \; u_x=g(x), \; u_y=f'(x)-g(x),
\nonumber\\
\hphantom{\gamma(x)=\bigg(}{}
u_{xx}=g'+\frac{g(f'-g)}{x-f}, \; u_{yy}=f''-g'+\frac{g(f'-g)}{x-f} \bigg),
\label{CHDE2}
\end{gather}
where $f(x)$ has no f\/ixed points.

The standard rank 3 Pfaf\/f\/ian system for the PDE in equation~\ref{pde1} is
given on a~7-mani\-fold~$M$ with coordinates $(x,y,u,u_x,u_y,u_{xx},u_{yy})$ by
\begin{gather*}
I=\spn\bigg\{\theta=du-u_xdx-u_ydy, \; \theta_x=du_x-u_{xx}dx-\frac{u_xu_y}{u-x}dy,
\\
\phantom{I=\spn\bigg\{ }
\theta_y=du_y-\frac{u_xu_y}{u-x}dx-u_{yy}dy\bigg\}.
\end{gather*}
With
\begin{gather*}
\homega=dx,\qquad\hpi=u_x d\left(\frac{u_{xx}}{u_x}+\frac{1}{u-x}\right),\qquad
\comega=dy,\qquad\cpi=u_y d\left(\frac{u_{yy}}{u_y}\right),
\end{gather*}
we have
\begin{gather*}
d\theta\equiv0,\qquad
d\theta_x\equiv\homega\wedge\hpi,\qquad
d\theta_y\equiv\comega\wedge\cpi\mod I,
\end{gather*}
and
\begin{gather*}
\hV^\infty=\spn\left\{ dx,\;d\left(\frac{u_x}{u-x}\right),\;d\left(\frac{u_{xx}}{u_x}+\frac{1}{u-x}\right) \right\},\qquad
\cV^\infty=\spn\left\{ dy,\;d\left(\frac{u_{yy}}{u_y}\right) \right\}.
\end{gather*}
The canonical quotient representation for $I$ can be found be found using the algorithm
in~\cite{anderson-fels-vassiliou:2009a}.
We f\/ind $M_1=\{(y,w,w_{y},w_{yy}),w_y>0)$, $M_2=\{(x,v,v_x,v_{xx},v_{xxx}),v_x>0\}$, and
\begin{gather}
K_1=\spn\{dw-w_y dy,\; dw_y-w_{yy}dy\},\nonumber
\\
K_2=\spn\{dv-v_x dx,\;dv_x-v_{xx}dx,\;dv_{xx}-v_{xxx}dx\}.
\label{KSE1}
\end{gather}
The Vessiot group $G$ is the 2-dimensional non-Abelian group $G=\{(a,b),\;a\in\Real^+,\; b\in\Real\}$
which acts on $M_1$ and~$M_2$ by
\begin{gather*}
\begin{split}
& (a,b)\cdot(y,w,w_y,w_{yy})=(y,aw-b,aw_y,a w_{yy}),
\\
& (a,b)\cdot(x,v,v_x,v_{xx},v_{xxx})=(x,av+b,av_x,av_{xx},a v_{xxx}).
\end{split}
\end{gather*}
The quotient map $\bfq_{G_\diag}:M_1\times M_2\to M$ can be written in coordinates as
\begin{gather}
\bfq(y,w,w_y, w_{yy};   x,v,v_x,v_{xx},v_{xxx} )
=\bigg(x=x,\;y=y,\;u=x-\frac{v+w}{v_x},
\nonumber\\
\qquad{} u_x=\frac{(v+w)v_{xx}}{v_x^2},\;u_y=-\frac{w_y}{v_x},\:u_{xx}=D_x(u_x),\;u_{yy}=\frac{-w_{yy}}{v_x}\bigg).
\label{CDE2a}
\end{gather}
A simple calculation shows that $\bfq_{G_\diag}^*(I)\subset K_1+K_2$ and~hence
$I=(K_1+K_2)/G_{\diag}$.
In this example this amounts to checking that
\begin{gather*}
u=x-\frac{v+w}{v_x}
\end{gather*}
from equation~\eqref{CDE2a} solves the PDE~\eqref{pde1}.

We proceed to solve the initial value problem using Theorem~\ref{slnthm}.
First we must f\/ind
the lift of the Cauchy data $\gamma$ given in~\eqref{CHDE2} to an integral curve $\sigma$ of $K_1+K_2$.
By using the second method of Remark~\ref{solvele}.
We
f\/ind the $3$-dimensional manifold $P=\bfq_{G_\diag}^{-1}(\gamma(x))\subset M_1\times M_2$
using~\eqref{CDE2a} and~\eqref{CHDE2} in terms of the parameters $x$, $v$, $v_x$,
\begin{gather}
P=\big(
y=x,\; w=(x-f(x))v_x-v,\; w_y=(g(x)-f'(x))v_x,
\nonumber\\
\hphantom{P=\big( }{}
w_{yy}=\left(g'(x)-f''(x)+(g(x)-f'(x))G(x)\right)v_x;
\nonumber\\
\hphantom{P=\big( }{}
x,\; v,\; v_x,\; v_{xx}=G(x)v_x,\; v_{xxx}=\big(G(x)^2+G'(x)\big)v_{x}
\big),
\label{PE1}
\end{gather}
where
\begin{gather}
G(x)=\frac{g(x)}{x-f(x)}.
\label{defG}
\end{gather}
The restriction of the Pfaf\/f\/ian system from equation~\eqref{KSE1} to $P$ in~\eqref{PE1} is then
\begin{gather}
(K_1+K_2)|_P=\spn\{ dv-v_x dx,\;d v_x-G(x)v_x dx \}.
\label{K1K2P1}
\end{gather}

We choose the point (see Remark~\ref{solvele})
\begin{gather*}
x_0=\bigg( x=0,\; y=0,\; u=f(0),\; u_x=g(0),\; u_y=f'(0)-g(0),
\\
\hphantom{x_0=\bigg( }{}
u_{xx}=g'(0)-\frac{g(0)(f'(0)-g(0))}{f(0)},\; u_{yy}=
f''(0)-g'(0)-\frac{g(0)(f'(0)-g(0))}{f(0)}\bigg)\in S
\end{gather*}
and $(x_1,x_2)\in M_1\times M_2$ to be (see~\eqref{PE1})
\begin{gather}
(x_1,x_2)=\big(y=0,\; w=-f(0),\; w_y=g(0)-f'(0),\nonumber\\
\hphantom{(x_1,x_2)=\big(}{}
w_{yy}=g'(0)-f''(0)+(f'(0)-g(0))G(0);
\nonumber\\
\hphantom{(x_1,x_2)=\big(}{} x=0,\; v=0,\; v_x=1,\; v_{xx}=G(0)
\big).
\label{ICDE2}
\end{gather}
These points satisfy $\bfq_{G_\diag}(x_1,x_2)=x_0=\gamma(0)$.
We now f\/ind the integral curve $\sigma$ of $(K_1+K_2)|_P$ through $(x_1,x_2)$.
This involves solving the Lie equation from~\eqref{K1K2P1} on $P$ subject to the initial
data~\eqref{ICDE2}.
With the $x$-coordinate as the parameter we f\/ind
\begin{gather}
v_x=e^{\int_0^x G(t)dt},\qquad v=\int_0^x e^{\int_0^s G(t)dt}ds.
\label{abdata3}
\end{gather}
The explicit form for $\sigma:\Real\to M_1\times M_2$ is found by inserting
equation~\eqref{abdata3} in~\eqref{PE1}, giving
\begin{gather}
\sigma(x)=\bigg( y=x,\; w=(x-f(x))e^{\int_0^x G(t)dt}-\int_0^x e^{\int_0^s G(t)dt}ds,
\;
w_y=(g(x)-f'(x))e^{\int_0^x G(t)dt},
\nonumber\\
\hphantom{\sigma(x)=\bigg(}{}
w_{yy}=\left(g'(x)-f''(x)+(g(x)-f'(x))G(x)\right)e^{\int_0^x G(t)dt};
\nonumber\\
\hphantom{\sigma(x)=\bigg(}{}
x=x,\;v=\int_0^x e^{\int_0^s G(t)dt}ds,\;v_x=e^{\int_0^x G(t)dt},\;v_{xx}=G(x)e^{\int_0^x G(t)dt},
\nonumber\\
\hphantom{\sigma(x)=\bigg(}{}
v_{xxx}=\big(G(x)^2+G'(x)\big)e^{\int_0^x G(t)dt}
\bigg).
\label{sige1}
\end{gather}

The curve $\sigma_1=\pi_1\circ\sigma$ is then easily determined from equation~\eqref{sige1}, which
in terms of the parameter $y$ is
\begin{gather}
\sigma_1(y)=\bigg( y,\;w=(y-f(y))e^{\int_0^y G(t)dt}-\int_0^y e^{\int_0^s G(t)dt}ds,\;w_y=(g(y)-f'(y))e^{\int_0^y G(t)dt},
\nonumber\\
\hphantom{\sigma_1(y)=\bigg(}{}
w_{yy}=\left(g'(y)-f''(y)+(g(y)-f'(y))G(y)\right)e^{\int_0^y G(t)dt}\bigg).
\label{abdata1}
\end{gather}
The curve $\sigma$ in equation~\eqref{sige1} also projects to the curve $\sigma_2=\pi_2\circ\sigma$
\begin{gather}
\sigma_2(x)=\bigg( x=x,\;v=\int_0^x e^{\int_0^s G(t)dt}ds,\;v_x=e^{\int_0^x G(t)dt},\;v_{xx}=
G(x)e^{\int_0^x G(t)dt}, \nonumber
\\
\hphantom{\sigma_2(x)=\bigg(}{}
v_{xxx}=\big(G'(x)+G(x)^2\big)e^{\int_0^x G(t)dt} \bigg).
\label{abdata2}
\end{gather}

Now, according to equation~\eqref{defSig} in Theorem~\ref{slnthm}, the solution to the
PDE~\eqref{pde1} with Cauchy data~\eqref{CHDE2} is the image of the product of the curves in
equation~\eqref{abdata1} and~\eqref{abdata2} in $M_1\times M_2$ under the map $\bfq_{G_\diag}$ in
equation~\eqref{CDE2a}.
This gives
\begin{gather}
u =x-\frac{v+w}{v_x}
=x-\frac{\int_0^x e^{\int_0^s G(t)dt}ds-(f(y)-y)e^{\int_0^y G(t)dt}-\int_0^y e^{\int_0^s G(t)dt}ds
}{e^{\int_0^x G(t)dt}}
\nonumber\\
\hphantom{u}
=x-\frac{\int_y^x e^{\int_0^s G(t)dt}ds-(f(y)-y)e^{\int_0^y G(t)dt}
}{e^{\int_0^x G(t)dt}}
\nonumber\\
\hphantom{u}
=x+(f(y)-y)e^{\int_x^y G(t)dt}+e^{-\int_0^x G(t)dt}\left(\int_x^y e^{\int_0^s G(t)dt}ds\right),
\label{E1sln}
\end{gather}
where $G(t)$ is given in equation~\eqref{defG}.
It is easy to check that this solves the Cauchy problem~\eqref{CHDE2} for the PDE~\eqref{pde1}.
One may interpret formula~\eqref{E1sln} as the analogue of d'Alembert's formula (in null
coordinates) for equation~\eqref{pde1}.

\end{Example}
\begin{Example} In this next example we write the standard rank 3 Pfaf\/f\/ian system $I$ for the
non Monge--Ampere hyperbolic equation
\begin{gather}
3u_{xx}u_{yy}^3+1=0
\label{eqex2}
\end{gather}
on a~$7$-manifold $M$ with coordinates $(x,y,u,u_x,u_y,u_{xy},u_{yy})$ as
\begin{gather*}
I=\spn\left\{ du-u_xdx-u_y dy,\;du_x+\frac{1}{3u_{yy}^3}dx-u_{xy}dy,\; du_y-u_{xy}dx-u_{yy}dy \right\}.
\end{gather*}

The construction of the canonical quotient representation for $I$ are given
in~\cite{anderson-fels:2011a} and~we summarize the result here.
On the $5$-manifolds $M_1$ with coordinates $(t,w,v,v_t,v_{tt})$ and~$M_2$ with coordinates
$(s,q,p,p_s,p_{ss})$ let
\begin{gather*}
\begin{split}
& K_1=\spn\big\{dw-v_{tt}^2dt,\; dv-v_t dt,\; dv_t-v_{tt}dt\big\},
\\
& K_2=\spn\big\{dq-p_{ss}^2ds,\; dp-p_sds,\; dp_s-p_{ss}ds\big\}.
\end{split}
\end{gather*}
The action of the group $G=\Real^3$ is given by
\begin{gather*}
(a,b,c)\cdot(t,w,v,v_t,v_{tt})=(t,w+a,v+b+ct,v_t+c,v_{tt}),
\\
(a,b,c)\cdot(s,q,p,p_s,p_{ss})=(s,q-a,p-b+cs,p_s+c,p_{ss}),\qquad a,b,c\in\Real.
\end{gather*}
The quotient map $\bfq_{G_\diag}:M_1\times M_2\to M$ is then given in these coordinates by
\begin{gather}
x=-2\frac{v_{tt}+p_{ss}}{s+t},\qquad
y=\frac{1}{2}(v_{tt}-p_{ss})(t+s)+p_{s}-v_{t},
\nonumber
\\
u=-q-w+2\frac{tv_t+sp_s-p-v}{s+t}(v_{tt}+p_{ss})
\nonumber
\\
\hphantom{u=}{}
+\frac{1}{3}\left(
(2s-t)v_{tt}^2+(2t-s)p_{ss}^2-2(s+t)v_{tt}p_{ss}\right),
\nonumber
\\
u_x=p+v-tv_t-sp_s+\frac{s+t}{6}\left((2t-s)v_{tt}+(2s-t)p_{ss}\right),
\nonumber
\\
u_y=2\frac{sv_{tt}-tp_{ss}}{s+t},
\qquad
u_{xy}=\frac{1}{2}(t-s),\qquad u_{yy}=\frac{2}{s+t}.
\label{Super3}
\end{gather}
It is straightforward matter to check $I=(K_1+K_2)/G_{\diag}$.

Consider the Cauchy data $\gamma:\Real\to M$ given by
\begin{gather*}
\gamma(\epsilon)=\big(x=0,\;y=\epsilon,\;u=f(\epsilon),\; u_x=g(\epsilon),\; u_y=f'(\epsilon),\; u_{xy}=
g'(\epsilon),\; u_{yy}=f''(\epsilon)\big).
\end{gather*}
The four-dimensional manifold $P=\bfq_{G_{\diag}}^{-1}(\gamma(\epsilon))\subset M_1\times M_2$ can
be computed from equation~\eqref{Super3}.
Using the parameters $\epsilon$, $w$, $v$, $v_t$, we get
\begin{gather*}
P=\Bigg( t=h(\epsilon),\;w,\;v,\;v_t,\;v_{tt}=\frac{1}{2}f';s=k(\epsilon),\;
q=\frac{(f')^2}{2f''}-w-f,
\\
\hphantom{P=\Bigg(}{}
p=g-\frac{f'}{(f'')^2}+\epsilon k(\epsilon)+2\frac{v_t}{f''}-v,\;
p_s=v_t+\epsilon-\frac{f'}{f''},\;p_{ss}=-\frac{1}{2}f' \Bigg),
\end{gather*}
where
\begin{gather*}
h(\epsilon)=\frac{1}{f''}+g' \qquad {\rm and}\qquad k(\epsilon)=\frac{1}{f''}-g'.
\end{gather*}
The restriction of $(K_1+K_2)|_P$ is
\begin{gather*}
(K_1+K_2)|_P=\spn\left\{ dw-\frac{(f')^2}{4}h'd\epsilon,\;dv-v_t h'd\epsilon,\;
dv_t-\frac{f'}{2}h'd\epsilon \right\}.
\end{gather*}
Taking the point $p_0=(\epsilon=0,v=0,w=0,v_t=0)\in P$ we f\/ind the maximal integral manifold $J$
of $(K_1+K_2)|_P$ through $p_0$ to be
\begin{gather}
v_t(\epsilon)=\frac{1}{2}\int_0^\epsilon h'f'd\xi,\qquad v(\epsilon)=
\int_0^\epsilon\left(\int_0^\tau h'f'd\tau\right)h'd\xi,\qquad w=\frac{1}{4}\int_0^\epsilon(f')^2h'd\xi.
\label{lcex2}
\end{gather}
Here $\bfq_{G_\diag}(p_0)=(0,0,u=f(0),u_x=g(0),u_y=g'(0),u_{yy}=f''(0))\in S$.
From Theorem~\ref{slnthm}, the solution to the initial value problem for~\eqref{eqex2}
is then found using~\eqref{lcex2} and~\eqref{Super3} to be
\begin{gather*}
x=\frac{f'(\delta)-f'(\epsilon)}{k(\delta)+h(\epsilon)},
\\
y=\delta+\frac{1}{4}f'(\epsilon)(h(\epsilon)+k(\delta))
+\frac{1}{4}f'(\delta)(h(\epsilon)-k(\delta))-\frac{1}{2}f'(\delta)h(\delta)
+\frac{1}{2}\int_\epsilon^\delta f'h'd\xi,
\\
u=f(\delta)-\frac{f'(\delta)^2}{2f''(\delta)}+\frac{1}{4}\int_{\epsilon}^{\delta}(f')^2h'd\xi+\frac{f'(\epsilon)-f'(\delta)}{h(\epsilon)+k(\delta)}
\left(\frac{f'(\delta)g'(\delta)}{f''(\delta)}-g(\delta)-\frac{1}{2}\int_\epsilon^\delta f'h'h d\xi\right)
\\
\hphantom{u=}{}
+\frac{k(\delta)}{12}(2f'(\epsilon)f'(\delta)+2f'(\epsilon)^2-f'(\delta)^2)+\frac{h(\epsilon)}{12}(2f'(\delta)^2-f'(\epsilon)^2+2f'(\epsilon)f'(\delta)).
\end{gather*}
\end{Example}

\begin{Example}
\label{lastE}
In this last example we use the system of two partial dif\/ferential equations
\begin{gather*}
u_{xz}=0,\qquad u_{yz}=0
\end{gather*}
to demonstrate the results in Section~\ref{decomp}.
The standard rank 4 Pfaf\/f\/ian system $I$ for this system on the $11$-manifold $M$ with
coordinates $(x,y,z,u,u_x,u_y,u_z,u_{xx},u_{xy},u_{yy},u_{zz})$ is given by
\begin{gather*}
I=\spn\big\{du-u_xdx-u_ydy-u_zdz,\;du_x-u_{xx}dx-u_{xy}dy,\\
\hphantom{I=\spn\big\{}{}
du_y-u_{xy}dx-u_{yy}dy,\; du_{z}-u_{zz}dz\big\}.
\end{gather*}
This is a~decomposable Pfaf\/f\/ian system (Section~\ref{decomp}) with $n_1=2$, $p_1=3$
and~$n_2=1$, $p_2=1$.
The canonical quotient representation for the Darboux integrable system~$I$ is given by taking
$M_1=J^2(\Real,\Real)$ and~$M_2=J^2(\Real^2,\Real)$ with
\begin{gather}
K_1 =\spn\{ dw-w_z dz,\;dw_z-w_{zz}dz \},
\nonumber\\
K_2 =\spn\{ dv-v_xdx-v_y dy,\;dv_x-v_{xx}dx-v_{xy}dy,\;dv_y-v_{xy}dx-v_{yy}dy \}.\label{K1K2E2}
\end{gather}
Then $I=(K_1+K_2)/G$ where the action of $G=\Real$ is given in coordinates by
\begin{gather*}
c\cdot(z,w,w_z,w_{zz})=(z,w+c,w_z,w_{zz}),
\\
c\cdot(x,y,v,v_x,v_y,v_{xx},v_{xy},v_{yy})=
(x,y,v-c,v_x,v_y,v_{xx},v_{xy},v_{yy}),\qquad c\in\Real.
\end{gather*}
The quotient map $\bfq_{G_\diag}:M_1\times M_2\to M$ written in the above coordinates is easily
found to be
\begin{gather}
\bfq_{G_{\diag}}(z,w,w_z,w_{zz};x,y,v,v_x,v_y,v_{xx},v_{xy},v_{yy})=
(x,\;y,\;z,\;u=w+v,\nonumber\\
\qquad{} u_x=v_x,\;u_y=v_y,\;u_z=w_z,\;u_{xx}=v_{xx},\;u_{xy}=v_{xy},\;u_{yy}=v_{yy},\;u_{zz}=w_{zz}).
\label{bfqE2}
\end{gather}

Let $a(x,y)$ be a~function of two variables and~$k(\xi)$ a~function of one, and~let $S$ be the
following non-characteristic two-dimensional integral manifold of $I$,
\begin{gather}
S=\big(x,\;y,\;z=x+y,\;u=a(x,y),\;u_x=a_x-k(x+y),\;u_y=a_y-k(x+y),
\nonumber\\
\hphantom{S=\big(}{}
u_z=k(x+y),\;u_{xx}=a_{xx}-k'(x+y),\;_{xy}=a_{xy}-k'(x+y),
\nonumber\\
\hphantom{S=\big(}{}
u_{yy}=a_{yy}-k'(x+y),\; u_{zz}=k'(x+y)\big).
\label{SE2}
\end{gather}
We proceed using Lemma~\ref{liftCC2}.
The set $P=\bfq_{G_\diag}^{-1}(S)$ is a~$3$-dimensional manifold.
With parameters $x$, $y$, $v$ it is easily determined from equations~\eqref{bfqE2} and~\eqref{SE2} to be
\begin{gather}
P=\big(z=x+y,\;w=a(x,y)-v,\;w_z=k(x+y),\;u_{zz}=k'(x+y);
\nonumber\\
\hphantom{P=\big(}{}
x,\;y,\;v,\;v_x=a_x-k(x+y),\;v_y=a_y-k(x+y),\;v_{xx}=a_{xx}-k'(x+y),
\nonumber\\
\hphantom{P=\big(}{}
v_{xy}=a_{xy}-k'(x+y),\;v_{yy}=a_{yy}-k'(x+y)\big).
\label{P3}
\end{gather}
On $P$ we have from~\eqref{K1K2E2},
\begin{gather}
(K_1+K_2)|_P=\spn\big\{dv-(a_x-k'(x+y))dx-(a_y-k'(x+y))dy\big\}.
\label{K1K2P2}
\end{gather}
We now f\/ind, as in Lemma~\ref{liftCC2}, the maximal $2$-dimensional integral manifold $L$ for the
Lie system in equation~\eqref{K1K2P2} subject to the initial conditions $(x=0,y=0,v(0,0)=a(0,0))$.
We obtain
\begin{gather}
v=a(x,y)-\int_0^{x+y}k(\xi)d\xi,
\label{VE2}
\end{gather}
and this determines $L\subset M_1\times M_2$ in Lemma~\ref{liftCC2}.
With $w=a(x,y)-v$ from equation~\eqref{P3}, where $v$ is given in~\eqref{VE2}, we f\/ind the
manifolds $L_i$ in Lemma~\eqref{fundl2} to be (with $n_1=2$ and~$n_2=1$)
\begin{gather}
L_1=\pi_1(L)=\bigg(z,\;w=\int_0^{z}k(\xi)d\xi,\;w_z=k(z),\;w_{zz}=k'(z)\bigg),
\nonumber\\
L_2=\pi_2(L)=\bigg(x,\;y,\;v=a(x,y)-\int_0^{x+y}k(\xi)d\xi,\;v_x=a_x-k(x+y),\;v_y=a_y-k(x+y),
\nonumber\\
\hphantom{L_2=\pi_2(L)=\bigg(}{}
 v_{xx}=a_{xx}\!-k'(x\!+\!y),\;v_{xy}=a_{xy}\!-k'(x\!+\!y),\;v_{yy}=a_{yy}\!-k'(x\!+\!y)\!\bigg).\!\!\!\!\!
\label{WE2}
\end{gather}
We now apply Theorem~\ref{slnthm3} to f\/ind the solution to the Cauchy problem.
With $u=w+v$ from equation~\eqref{bfqE2}
and substituting for $w$ and~$v$ in~\eqref{WE2} we have
\begin{gather*}
u(x,y,z)=\int_{0}^z k(\xi)d\xi+a(x,y)-\int_0^{x+y}k(\xi)d\xi
=a(x,y)+\int_{x+y}^z k(\xi)d\xi.
\end{gather*}
\end{Example}

\appendix

\section{The classical theory and~a~second approach}
\label{CCP}

In this appendix we recall the classical theory for solving the Cauchy problem for Darboux
integrable hyperbolic systems as given in Remark~9.4 of~\cite{anderson-fels:2011a} (demonstrated
in~\cite{bryant-griffiths-hsu:1995a}) and~then relate this approach to the approach taken to
solve the initial value problem given by Theorem~\ref{slnthm}.
This will lead to an alternative but equivalent method of solution to the Cauchy problem which
again requires solving equations of fundamental Lie type.
Example~\ref{Ex3} below demonstrates the theory.

We denote by $\CalI$ the Pfaf\/f\/ian exterior dif\/ferential system (EDS) generated by the
sections of~$I$.
If $\bfp:M\to N$ is a~surjective submersion we def\/ine the reduction of $\CalI$ by $\bfp$ as
\begin{gather}
\CalI/\bfp=\big\{ \theta\in\Omega^*(N)\,|\,\bfp^*\theta\in\CalI\big\}.
\label{Ip}
\end{gather}
If $G$ is a~symmetry group then $\CalI/G=\CalI/\bfq_G$.
For more details on EDS reduction see~\cite{anderson-fels:2011b}.

Let $\CalI$ be a~Darboux integrable hyperbolic Pfaf\/f\/ian system and~let $G$, $M_i$, $\CalK_i$ be the
data for the canonical quotient representation of $\CalI$ as in Theorem~\ref{Dreduce}
or~\ref{Intro4}.
Let $\pi_i:M_1\times M_2\to M_i$ and~let $\bfq_G^i:M_i\to M_i/G$ be the quotient maps.
The compositions $\bfq_G^1\circ\pi_1$
and $\bfq_G^2\circ\pi_2$ are invariant with respect to the diagonal action of $G$ on $M_1\times
M_2$ and~therefore these maps factor through $\bfq_{G_\diag}$.
Accordingly we can def\/ine the surjective submersions $\bfp_i:M\to M_i/G$ so that the following
diagram commutes (see also~\cite[equation~(6.6), p.~21]{anderson-fels:2011a}):
\begin{gather}
\begin{gathered}
\beginDC{\commdiag}[3]
\obj(0,20)[I]{$\left(\CalK_1+\CalK_2,M_1\times M_2\right)$}
\obj(0,0)[H]{$\left(\CalI,M\right)$}
\obj(50,20)[Ma]{$\left(\CalK_2,M_2\right)$}
\obj(50,0)[Maq]{$\left(\CalK_2/G,M_2/G\right).$}
\obj(-50,20)[Ma1]{$\left(\CalK_1,M_1\right)$}
\obj(-50,0)[Maq1]{$\left(\CalK_1/G,M_1/G\right)$}
\mor{I}{H}{$\bfq_{G_\diag}$}[\atleft,\solidarrow]
\mor{H}{Maq}{$\bfp_2$}[\atleft,\solidarrow]
\mor{I}{Ma}{$\pi_2$}[\atleft,\solidarrow]
\mor{Ma}{Maq}{$\bfq_G^2$}[\atleft,\solidarrow]
\mor{H}{Maq1}{$\bfp_1$}[\atright,\solidarrow]
\mor{I}{Ma1}{$\pi_1$}[\atright,\solidarrow]
\mor{Ma1}{Maq1}{$\bfq_G^1$}[\atleft,\solidarrow]
\enddc
\end{gathered}
\label{CDP}
\end{gather}

The commutativity of diagram~\eqref{CDP} means that the Pfaf\/f\/ian systems $\CalK_i$ satisfy (see
Theorem~3.1 in~\cite{anderson-fels:2011a})
\begin{gather}
\CalK_i/G=\CalI/\bfp_i\qquad{\rm and}\qquad(\CalK_1+\CalK_2)/\pi_i=\CalK_i.
\label{EDSC}
\end{gather}
It is important to note that the maps $\bfq_{G_{\diag}}:(\CalK_1+\CalK_2,M_1\times
M_2)\to(\CalI,M)$ and~$\bfq^i_G:(\CalK_i,M_i)\to(\CalK_i/G,M_i/G)$ in diagram~\eqref{CDP} all
satisfy the hypothesis of Theorem~\ref{FRD} except that $G$ may not act transversally to $K_i'$.
In this case~$\CalK_i/G$ will not be a~Pfaf\/f\/ian system.
This is demonstrated in Example~\ref{Ex3}.
The maps $\bfp_i:(\CalI,M)\to(\CalK_i/G,M/G_i)$ are shown in~\cite{anderson-fels:2011a} to be
integrable extensions~\cite{bryant-griffiths:1995p}.
In particular the maps $\bfp_i$ map integral manifolds of $\CalI$ to integral manifolds of~$\CalK_i/G$.

\subsection{The classical method}
\label{Classm}

To describe the classical integration method we assume that the two integrable distribu\-tions~$\cV^\infty$ and~$\hV^\infty$ are regular so that the two maps to the leaf spaces
$\tilde\bfp_1:M\to M/\ann(\cV^\infty)$ and
$\tilde\bfp_2:M\to M/\ann(\hV^\infty)$ are smooth submersions.
By using def\/inition~\eqref{Ip} we can then construct the following diagram
\begin{gather}
\beginDC{\commdiag}[3]
\obj(-45,0)[I]{$\big(\CalI/\tilde\bfp_1,M/\ann\big(\cV^\infty\big)\big)$}
\obj(0,0)[H]{$\left(\CalI,M\right)$}
\obj(45,0)[Ma]{$\big(\CalI/\tilde\bfp_2,M/\ann\big(\hV^\infty\big)\big).$}
\mor{H}{I}{$\tilde\bfp_1$}[\atright,\solidarrow]
\mor{H}{Ma}{$\tilde\bfp_2$}[\atleft,\solidarrow]
\enddc
\label{subdiag}
\end{gather}
The quotient spaces $M/\ann(\hV^\infty)$ and~$M/\ann(\cV^\infty)$ are thought of as the spaces of
intermediate integrals in the sense that
\begin{gather*}
\cV^\infty=\tilde\bfp^*_1\big(T^*(M/\ann(\cV^\infty))\big)\qquad{\rm and}\qquad\hV^\infty=
\tilde\bfp^*_2\big(T^*(M/\ann(\hV^\infty))\big).
\end{gather*}
The f\/irst of these equations implies that if $F:M/\ann(\cV^\infty)\to\Real$, then
$\tilde\bfp_1^*(dF)$ takes values in~$\cV^\infty$, and~that a~local basis of sections for
$\cV^\infty$ can be constructed in this way.
A similar statement holds for $\hV^\infty$.

Let $\gamma:(a,b)\to M$ be a~non-characteristic integral curve of the Darboux integrable hyperbolic
system $\CalI$.
We project the curve $\gamma$ into the spaces of intermediate integrals by
\begin{gather}
\tilde\gamma_i=\tilde\bfp_i\circ\gamma.
\label{deflamit}
\end{gather}
The $\gamma_i$ are integral curves of $\CalI/\tilde\bfp_i$.
Now let $N\subset M$ be the inverse image of the product curves
\begin{gather*}
N=(\tilde\bfp_1,\tilde\bfp_2)^{-1}\left(\tilde\gamma_1(a,b),\tilde\gamma_2(a,b)\right).
\end{gather*}
Then $\CalI|_{N}$ is a~Frobenius system.
By determining the maximal integral manifold of $\CalI|_{N}$ through a~point $\gamma(t_0)$ of
$\CalI|_{N}$, one determines a~local solution to the initial value problem.
This is the classical method.

It is important to note that in the classical method the Vessiot group $G$ does not appear
and~there is no reason to expect that the Frobenius system $\CalI|_{N}$ can be integrated by using
equations of Lie type.

\subsection{A second method and~the comparison}

The relationship between the classical method and~Theorem~\ref{slnthm} can be described using
diagram~\eqref{CDP} once we make one key observation.
Namely that {\it diagram~\eqref{subdiag} can be identified with the bottom row of
diagram~\eqref{CDP}.} To show how that identif\/ication can be made, suppose that the manifolds
$M_i$ are connected, then the spaces of intermediate integrals $M/\ann(\cV^\infty)$
and~$M/\ann(\hV^\infty)$ may be identif\/ied with the two quotient spaces $M_i/G$.
In particular it is shown in Theorem~6.1 of~\cite{anderson-fels:2011a} that
\begin{gather}
\cV^\infty=\bfp_1^*\big(T^*(M_1/G)\big)\qquad{\rm and}\qquad
\hV^\infty=\bfp_2^*\big(T^*(M_2/G)\big).
\label{CIV}
\end{gather}
From diagram~\eqref{CDP}, part $(v)$ of Theorem~\ref{Dreduce}, and~\eqref{redannI} we have
\begin{gather*}
\ker(\bfp_{1*})=\bfq_{G_\diag*}(0+TM_2)=\ann\big(\cV^\infty\big).
\end{gather*}
The f\/ibres of $\bfp_1$ are of the form $\bfq_{G_\diag}(x_1,M_2)$ which are connected since $M_2$
is connected.
Therefore by Theorem 3.18 in~\cite{anderson-fels:2011b},
\begin{gather}
M/\ann\big(\cV^\infty\big)\cong M_1/G\qquad{\rm and}\qquad\CalI/\tilde\bfp_1\cong\CalK_1/G,
\label{cong1}
\end{gather}
where the f\/irst equivalence is by canonical dif\/feomorphism, and~the second equivalence follows
from the f\/irst equation in~\eqref{EDSC}.
Similarly we have
\begin{gather}
M/\ann\big(\hV^\infty\big)\cong M_2/G\qquad{\rm and}\qquad\CalI/\tilde\bfp_2\cong\CalK_2/G.
\label{cong2}
\end{gather}

With the identif\/ications~\eqref{cong1} and~\eqref{cong2}, we continue using diagram~\eqref{CDP}
to relate the classical solution to the Cauchy problem described above and~Theorem~\ref{slnthm}.
From the initial data $\gamma$ def\/ine the maps $\gamma_i:(a,b)\to M_i/G$ by
\begin{gather}
\gamma_i=\bfp_i\circ\gamma.
\label{deflami}
\end{gather}
Note that with the identif\/ications in equation~\eqref{cong1} and~\eqref{cong2}, the curves
$\tilde\gamma_i$ in~\eqref{deflamit} and~$\gamma_i$ in~\eqref{deflami} are identif\/ied.
The following lemma follows from the fact that the maps $\bfp_i:(\CalI,M)\to(\CalK_i/G,M_i)$ are
integrable extensions.
\begin{Lemma}
The curves $\gamma_i:(a,b)\to M_i/G$ are integral curves of $\CalK_i/G$.
\end{Lemma}

Recall now that in Theorem~\ref{slnthm} we choose a~point $x_0=\gamma(t_0)$ and~a~point
$(x_1,x_2)\in M_1\times M_2$ with $\bfq_{G_\diag}(x_1,x_2)=x_0$.
We then let $\sigma:(a,b)\to M_1\times M_2$ be the lift of $\gamma$ obtained by Theorem~\ref{FRD}
through the point $\sigma(t_0)=(x_1,x_2)$.
With
\begin{gather*}
\sigma_i=\pi_i\circ\sigma,
\end{gather*}
def\/ined in equation~\eqref{sigmai} in Lemma~\ref{fundl},
the solution to the initial value problem is given by
$s(t_1,t_2)=\bfq_{G_\diag}\big(\sigma_1(t_1),\sigma_2(t_2)\big)$ in equation~\eqref{defSig}.

Putting the curves $\gamma$, $\sigma$, $\sigma_i$ and~$\gamma_i$ into the commutative
diagram~\eqref{CDP}, we have
\begin{gather}
\beginDC{\commdiag}[3]
\obj(-35,10)[I]{$\left(\CalK_1+\CalK_2,M_1\times M_2\right)$}
\obj(-35,-10)[H]{$\left(\CalI,M\right)$}
\obj(35,10)[Ma]{$\left(\CalK_i,M_i\right)$}
\obj(35,-10)[Maq]{$\left(\CalK_i/G,M_i/G\right),$}
\obj(0,0)[ab]{$(a,b)$}
\mor{ab}{H}{$\gamma$}[\atleft,\solidarrow]
\mor{ab}(-33,8){$\sigma$}[\atleft,\solidarrow]
\mor{ab}{Ma}{$\sigma_i$}[\atright,\solidarrow]
\mor{ab}(33,-8){$\gamma_i$}[\atright,\solidarrow]
\mor{H}{Maq}{$\bfp_i$}[\atright,\solidarrow]
\mor{I}{Ma}{$\pi_i$}[\atright,\solidarrow]
\mor{Ma}{Maq}{$\bfq_G^i$}[\atleft,\solidarrow]
\mor{I}{H}{$\bfq_{G_\diag}$}[\atright,\solidarrow]
\enddc
\label{CDP2}
\end{gather}
Since $\bfq_G^i\circ\sigma_i=\gamma_i$ and~the initial points of $\sigma_i$ and~$\gamma_i$
satisfy
\begin{gather}
\sigma_i(t_0)=x_i,\qquad\gamma_i(t_0)=\bfp_i(\bfq_{G_\diag}(x_1,x_2))=\bfp_i(x_0)=
\bfq_G^i(\sigma_i(t_0)),
\label{ICS}
\end{gather}
we have by Theorem~\ref{FRD} and~the commutative diagram~\eqref{CDP2} that $\sigma_i$ is the unique
lift of $\gamma_i$ through $x_i$ which is an integral curve of $\CalK_i$.
{\it Consequently the curves $\sigma_i$ in Theorem~{\rm \ref{slnthm}} can be found from the curves
$\gamma_i$, which are determined from the initial data $\gamma$ in~\eqref{deflami}, by solving
two separate equations of fundamental Lie type.} This gives us the following alternate method of
solving the initial value problem.
\begin{Theorem}
\label{slnthm2}
Let $\gamma:(a,b)\to M$ be a~non-characteristic integral curve for the Darboux integrable,
hyperbolic Pfaffian system $I$ with Vessiot group $G$.
Let $\gamma_i:(a,b)\to M_i/G$ be the projection of $\gamma$ to the two spaces~$M_i/G$, and~let
$\sigma_i:(a,b)\to M_i$ be the unique lift of~$\gamma_i$ to an integral curve of $\CalK_i$ on $M_i$
satisfying the initial conditions in~\eqref{ICS}.
Each curve $\sigma_i$ is found by solving an equation of fundamental Lie type and~the function
\begin{gather}
s(t_1,t_2)=\bfq_{G_\diag}\big(\sigma_1(t_1),\sigma_2(t_2)\big),\qquad t_1,t_2\in(a,b)
\label{defSig2}
\end{gather}
solves the Cauchy problem for~$\gamma$.
\end{Theorem}

At this point we can compare the approach to solving the Cauchy problem given in
Theorem~\ref{slnthm2} to the classical method in Section~\ref{Classm}.
Starting with the Cauchy data $\gamma:(a,b)\to M$ to~$\CalI$
we project~$\gamma$ to curves $\tilde\gamma_i:(a,b)\to M_i/G$.
The classical method then constructs the solution to the initial value problem by integrating the
Frobenius system $\CalI|_{N}$ where
\begin{gather*}
N=(\bfp_1,\bfp_2)^{-1}\left(\gamma_1(a,b),\gamma_2(a,b)\right).
\end{gather*}
On the other hand Theorem~\ref{slnthm2} constructs the solution by using the quotient
representation of~$\CalI$ and~lifting the projected curves $\gamma_i$ to integral curves of the EDS
$\CalK_i$ on~$M_i$.
The lifting problem involves solving equations of fundamental Lie type for the Vessiot group~$G$.
It is at this point we see the fundamental dif\/ference between the classical method of solving the
Cauchy problem and~the one introduced here in Theorem~\ref{slnthm2}.
This comparison shows how the Frobenius system in the classical solution can be integrated using
equations of Lie type.
As a~consequence of our approach we can identify when the Cauchy problem for generic initial data
can be solved by quadratures.

It is also worth noting that the projected curves~$\gamma_i$ in equation~\eqref{deflami} in the
classical approach are thought of (through equation~\eqref{CIV}) as prescribing the value of the
intermediate integrals for $\CalI$ along the initial data.
Integrating $\CalI|_N$ in Section~\ref{Classm} is then considered to be solving the ``prescribed
intermediate integral'' problem.
In Theorem~\ref{slnthm2} we view the projected initial data $\gamma_i$ as prescribing ``curvature
invariants'' of the curves $\sigma_i$ which we are trying to f\/ind.
See Example~\ref{Ex3} below where this is demonstrated explicitly.

As our f\/inal remark we compare the two solutions to the Cauchy problem in Theorems~\ref{slnthm}
and~\ref{slnthm2}.
In Theorem~\ref{slnthm} we lift the Cauchy data $\gamma$ to $M_1\times M_2$ and~then projects it to
curves $\sigma_i:(a,b)\to M_i$.
The solution $s$ is then constructed from these curves $\sigma_i$ in~\eqref{defSig}.
In the approach Theorem~\ref{slnthm2}, we project the Cauchy data~$\gamma$ to curves
$\gamma_i:(a,b)\to M_i/G$ and~then lift these curves to $\tilde\sigma_i:(a,b)\to M_i$.
It is not dif\/f\/icult to show that the corresponding solutions $s$ in~\eqref{defSig2}
and~\eqref{defSig} are then identical.

\begin{Example}
\label{Ex3}
In this example we write the standard rank 3 Pfaf\/f\/ian system $I$ for the Liouville equation
\begin{gather*}
u_{xy}=e^u
\end{gather*}
on a~$7$-manifold $M$ with coordinates $(x,y,u,u_x,u_y,u_{xy},u_{yy})$ as
\begin{gather*}
I=\spn\big\{ du-u_xdx-u_y dy,\;du_x-u_{xx}dx-e^u dy,\; du_y-e^u dx-u_{yy}dy\big\}.
\end{gather*}
The intermediate integrals are given by
\begin{gather*}
\cV^\infty=\spn\left\{ dy,\;d\left(u_{yy}-\frac{1}{2}u_y^2\right)\right\}\qquad{\rm and}\qquad
\hV^\infty=\spn\left\{ dx,\;d\left(u_{xx}-\frac{1}{2}u_x^2\right)\right\}.
\end{gather*}

The details of the canonical quotient representation for $I$ are given
in~\cite{anderson-fels:2011a} and~we summarize them here.
Let $K_1$ and~$K_2$ be the standard contact system on $J^3(\Real,\Real)$ and~$J^3(\Real,\Real)$.
In local coordinates $(y,w,w_y,w_{yy},w_{yyy})$ and~$(x,v,v_x,v_{xx},v_{xxx})$ we have
\begin{gather*}
K_1 =\spn\big\{ dw-w_ydy,\;d w_y-w_{yy}dy,\;dw_{yy}-w_{yyy}dy\big\}\qquad{\rm and}
\\
K_2 =\spn\big\{ dv-v_x dx,\;d v_x-v_{xx}dx,\;dv_{xx}-v_{xxx}dx\big\}.
\end{gather*}
Let
\begin{gather*}
\Gamma_1=\spn\big\{ \partial_w,\;\pr(w\partial_w),\;\pr\big(w^2\partial_{w}\big)\big\}\qquad{\rm and}\qquad
\Gamma_2=\spn\big\{ \partial_v,\;\pr(v\partial_v),\;\pr\big(v^2\partial_{v}\big)\big\}
\end{gather*}
be the Lie algebras of vector-f\/ields given by the prolongation of the standard inf\/initesimal
action of $\sl(2,\Real)$ acting on $w$ and~$v$, and~let
\begin{gather*}
\Gamma_{\diag}=
\big\{\partial_w-\partial_v,\;\pr(w\partial_w)+\pr(v\partial_v),\;\pr\big(w^2\partial_w\big)-\pr(v^2\partial_v)\big\}
\end{gather*}
denote the diagonal action.
On $M$, the open set where $v,v_x,w,w_y>0$, the distribution $\Gamma_\diag$ is regular
and~$M/\Gamma_{\diag}$ is 7-dimensional.
The quotient map $\bfq_{\Gamma_{\diag}}$ can be written in coordinates as
\begin{gather}
\bfq_{\Gamma_{\diag}}=\bigg( x=x,\;y=y,\;u=\log\frac{2w_y v_x}{v+w},\;u_x=
\frac{v_{xx}}{v_x}-2\frac{v_x}{v+w},
\nonumber\\
\hphantom{\bfq_{\Gamma_{\diag}}=\bigg(}{}
u_y=\frac{w_{yy}}{w_y}-2\frac{w_y}{v+w},
\; u_{xx}=\frac{v_{xxx}}{v_x}-\frac{v_{xx}^2}{v_x^2}-\frac{2v_{xx}}{v+w}+\frac{2v_x^2}{(v+w)^2},
\nonumber\\
\hphantom{\bfq_{\Gamma_{\diag}}=\bigg(}{}
u_{yy}=\frac{w_{yyy}}{w_y}-\frac{w_{yy}^2}{w_y^2}-\frac{2w_{yy}}{v+w}+\frac{2w_y^2}{(v+w)^2}
\bigg).
\label{bfp1}
\end{gather}

We now construct the lower part of diagram~\eqref{CDP} by computing the quotients
$(\CalK_i/{\Gamma_i},M_i/\Gamma_i)$.
The projection maps into the dif\/ferential invariants of $\Gamma_i$ on $M_i$ are
\begin{gather}
\bfq_{\Gamma_1}(y,w,w_y,w_{yy},w_{yyy}) =\left( y=y,\;\widecheck s=
\frac{w_{yyy}}{w_y}-\frac{3w_{yy}^2}{2w_y^2} \right)\qquad{\rm and}
\nonumber\\
\bfq_{\Gamma_2}(x,v,v_x,v_{xx},v_{xxx}) =\left( x=x,\;\widehat s=
\frac{v_{xxx}}{v_x}-\frac{3v_{xx}^2}{2v_x^2} \right),
\label{MmodGamma}
\end{gather}
where $\widecheck s$ and~$\widehat s$ are the Schwarzian derivatives of $w(y)$ and~$v(x)$.
In terms of the dif\/ferential invariants~$\widecheck s$ and~$\widehat s$, the algebraic generators
of the dif\/ferential systems~$\CalK_i$ can be written
\begin{gather}
\CalK_1 =
\langle dw-w_y dy,\;d w_y-w_{yy}dy,\;dw_{yy}-w_{yyy}dy, \; d\widecheck s\wedge dy\rangle_{\text{alg}},
\nonumber\\
\CalK_2 =
\langle dv-v_x dx,\;d v_x-v_{xx}dx,\;dv_{xx}-v_{xxx}dx,\; d\widehat s\wedge dx\rangle_{\text{alg}}.
\label{MFK1}
\end{gather}
The quotients are then quickly computed from~\eqref{MFK1} in the coordinates from~\eqref{MmodGamma}
to be
\begin{gather}
\CalK_1/\Gamma_1=\langle d\widecheck s\wedge dy\rangle \qquad{\rm and}\qquad
\CalK_2/\Gamma_2=\langle d\widehat s\wedge dx\rangle.
\label{KmodGam}
\end{gather}
Equation~\eqref{KmodGam} shows that an integral curve of $\CalK_1/\Gamma_1$ can be simply thought
of as a~choice of projective curvature $\widecheck s=G(y)$ and~an integral curve of
$\CalK_2/{\Gamma_2}$ is another choice of projective curvature $\widehat s=F(x)$.

The projection maps $\bfp_i:M\to M_i/\Gamma_i$ can be determined using the coordinates from
equations~\eqref{bfp1} and~\eqref{MmodGamma}.
We f\/ind
\begin{gather}
\widecheck s=\frac{w_{yyy}}{w_y}-\frac{3w_{yy}^2}{2w_y^2}=u_{yy}-\frac{1}{2}u_{y}^2,
\qquad
\widehat s =\frac{v_{xxx}}{v_x}-\frac{3v_{xx}^2}{2v_x^2}=u_{xx}-\frac{1}{2}u_{x}^2,
\label{iie}
\end{gather}
and so
\begin{gather}
\bfp_1(x,y,u,u_{x},u_{y},u_{xx},u_{yy}) =\left(y=y,\; \widecheck s=
u_{yy}-\frac{1}{2}u_{y}^2\right),
\nonumber\\
\bfp_2(x,y,u,u_{x},u_{y},u_{xx},u_{yy}) =\left(x=x,\; \widehat s=u_{xx}-\frac{1}{2}u_{x}^2\right).
\label{rmp}
\end{gather}
Equations~\eqref{iie} and~\eqref{rmp} shows that the (space of) intermediate integrals~$\cV^\infty$
and~$\hV^\infty$ are the pullback of the corresponding (space of) dif\/ferential invariants of~$\Gamma_i$ on~$M_i$.
This is the content of equation~\eqref{CIV}.

We now consider the initial value problem for the non-characteristic integral manifold
$\gamma:\Real\to M$ of $I$ given by
\begin{gather}
\gamma(x)=\big( x=x,\;y=x,\;u=f(x),\;u_x=g(x),\;u_y=f'(x)+g(x),\nonumber
\\
\hphantom{\gamma(x)=\big(}{}
u_{xx}=g'(x)-e^{f(x)},\; u_{yy}=f''(x)+g'(x)-e^{f(x)}\big).
\label{gamdata3}
\end{gather}
The projected curves $\gamma_i=\bfp_i\circ\gamma:M\to M_i/\Gamma_i$ are determined from~\eqref{rmp}
and~\eqref{gamdata3} to be
\begin{gather}
\gamma_1(y)=( y,\;\widecheck s=G(y) )\qquad{\rm and}\qquad
\gamma_2(x)=( x,\;\widehat s=F(x) ),
\label{deflami3}
\end{gather}
where
\begin{gather}
G(y)=f''(y)+g'(y)-e^{f(y)}-\frac{1}{2}(f'(y)+g(y))^2\qquad{\rm and}\nonumber\\
F(x)=g'(x)-e^{f(x)}-\frac{1}{2}g(x)^2.
\label{FandG}
\end{gather}

The classical method of solving the initial value problem~\eqref{gamdata3} given in
Section~\ref{Classm} is to f\/ind the integral manifolds of the system of~$\CalI$ restricted to
$(\bfp_1,\bfp_2)^{-1}(y,G(y),x,F(x))$.
This restricted system is always completely integrable.
By equations~\eqref{deflami},~\eqref{rmp} and~\eqref{FandG} we have
\begin{gather*}
(\bfp_1,\bfp_2)^{-1}(y,G(y),x,F(x))=
\left(x,y,u,u_{x},u_{y},u_{xx}-\frac{1}{2}u_{x}^2=F(x),u_{yy}-\frac{1}{2}u_{y}^2=G(y)\right),
\end{gather*}
where $G(y)$ and~$F(x)$ are given by~\eqref{FandG}.
The restriction of $\CalI$ to this subset is the completely integrable Pfaf\/f\/ian system
corresponding to the over-determined system of PDE
\begin{gather}
u_{xx}-\frac{1}{2}u_x^2 =F(x)=g'(x)-e^{f(x)}-\frac{1}{2}g(x)^2,
\qquad
u_{xy} =e^u,
\nonumber\\
u_{yy}-\frac{1}{2}u_{y}^2 =G(y)=f''(y)+g'(y)-e^{f(y)}-\frac{1}{2}(f'(y)+g(y))^2.
\label{Frob3}
\end{gather}
The solution to~\eqref{Frob3} solves the initial value problem~\eqref{gamdata3} for the Darboux
integrable equation $u_{xy}=e^u$.

The alternative method for solving the initial value problem~\eqref{gamdata3}
is presented in Theorem~\ref{slnthm2}.
In this case we must construct integral curves $\sigma_i:\Real\to M_i$ of the contact systems~$\CalK_i$ such that $\bfq_{\Gamma_i}\circ\sigma_i=\gamma_i$, where $\gamma_i$ are given
in~\eqref{deflami3}.
Using the expression for~$\widecheck s$ and~$\widehat s$ from equation~\eqref{MmodGamma} we seek
curves $w(y)$ and~$v(x)$ such that
\begin{gather}
\frac{w_{yyy}}{w_y}-\frac{3w_{yy}^2}{2w_y^2} =G(y)=f''(y)+g'(y)-e^{f(y)}-\frac{1}{2}(f'(y)+g(y))^2,
\nonumber\\
\frac{v_{xxx}}{v_x}-\frac{3v_{xx}^2}{2v_x^2} =F(x)=g'(x)-e^{f(x)}-\frac{1}{2}g(x)^2.
\label{Rconstruct}
\end{gather}
A solution to~\eqref{Rconstruct} produces a~solution to the initial value problem~\eqref{gamdata3}
by determining the curves~$\sigma_i$ in the formula~\eqref{defSig2} in Theorem~\ref{slnthm2}.

It is worth pointing out that solving equations~\eqref{Rconstruct} is a~classical reconstruction
problem in projective dif\/ferential geometry.
That is, given the Schwarzian or projective curvature, $\widehat s=F(x)$, $x\in(c,d)$, construct
a~curve $\alpha:(c,d)\to\Real{\rm P}^1$ such that the Schwarzian derivative of~$\alpha$ is~$F(x)$.
This is well known to be equivalent to solving an equation of fundamental Lie type for the simple
Lie group ${\rm PSL}(2,\Real)$.
\end{Example}

\pdfbookmark[1]{References}{ref}
\LastPageEnding

\end{document}